\documentclass[10pt, a4paper, reqno]{amsart}

\newtheorem{Theorem}{Theorem}
\newtheorem{Proposition}{Proposition}
\newtheorem{Lemma}{Lemma}
\newtheorem{Corollary}{Corollary}
\newtheorem{Remark}{Remark}
\newtheorem{Example}{Example}
\newtheorem{Definition}{Definition}

\author{Nevio Dubbini \and Benedetto Piccoli \and Antonio Bicchi}
\address{Nevio Dubbini: Interdepartmental research center ``E.
Piaggio'', University of Pisa, Italy, (e-mail:
nevio.dubbini@for.unipi.it)}

\address{Benedetto Piccoli: Department of mathematical sciences,
Rutgers University, New Jersey, USA}

\address{Antonio Bicchi: Interdepartmental center ``E. Piaggio'',
University of Pisa, Italy}

\date{}

\usepackage{graphics} 
\usepackage{epsfig} 
\usepackage{amsmath} 
\usepackage{amssymb}  

\title[Left invertibility of discrete-time output-quantized systems]{Left invertibility of discrete-time output-quantized systems: the linear case with finite inputs}

\begin{document}

\maketitle

\begin{abstract}
This paper studies left invertibility of discrete-time linear
output-quantized systems. Quantized outputs are generated according
to a given partition of the state-space, while inputs are sequences
on a finite alphabet. Left invertibility, i.e. injectivity of I/O
map, is reduced to left D-invertibility, under suitable conditions.
While left invertibility takes into account membership to sets of a
given partition, left D-invertibility considers only membership to a
single set, and is much easier to detect. The condition under which
left invertibility and left D-invertibility are equivalent is that
the elements of the dynamic matrix of the system form an
algebraically independent set. Our main result is a method to
compute left D-invertibility for all linear systems with no
eigenvalue of modulus one. Therefore we are able to check left
invertibility of output-quantized linear systems for a full measure
set of matrices. Some examples are presented to show the application
of the proposed method.
\end{abstract}

\section{Introduction}

Left invertibility is an important problem of systems theory, which
corresponds to injectivity of I/O map. It deals so with the
possibility of recovering unknown inputs applied to the system from
the knowledge of the outputs.

We investigate left invertibility of discrete--time linear
output-quantized systems in a continuous state-space. In particular,
inputs are arbitrary sequences of symbols in a finite alphabet: each
symbol is associated to an action on the system. Information
available on the system is represented by sequences of output
values, generated by the system evolution according to a given
partition of the state-space (quantization).

In recent years there has been a considerable amount of work on
quantized control systems (see for instance \cite{Del,Pica,Szn} and
references therein), stimulated also by the growing number of
applications involving ``networked'' control systems, interconnected
through channels of limited capacity (see e.g.
\cite{Bic,Car,Tat,Van}). The quantization and the finite cardinality
of the input set occur in many communication and control systems.
Finite inputs arise because of the intrinsic nature of the actuator,
or in presence of a logical supervisor, while output quantization
may occur because of the digital nature of the sensor, or if data
need a digital transmission.

Applications of left invertibility include fault detection in
Supervisory Control and Data Acquisition (SCADA) systems, system
identification, and cryptography (\cite{Edel,Inou}). Invertibility
of linear systems is a well understood problem, first handled in
\cite{Brok}, and then considered with algebraic approaches (see e.g.
\cite{Silv}), frequency domain techniques (\cite{Mas,Mas2}), and
geometric tools (cf. \cite{Mors}). Invertibility of nonlinear
systems is discussed in \cite{Resp}. More recent work has addressed
the left invertibility for switched systems (\cite{Vu}), and for
quantized contractive systems (\cite{Nev}).

The main intent of the paper is to show that the analysis of left
invertibility can be substituted, under suitable conditions, by an
analysis of a stronger notion, called left D-invertibility. The
condition under which left invertibility and left D-invertibility
are equivalent is that the elements of the dynamic matrix of the
system form an algebraically independent set (Theorem \ref{Theor alg
independent}). Therefore the set of matrices for which left
D-invertibility and left invertibility are equivalent is a full
measure set. While left invertibility takes in account whether two
states are in the same element of a given partition, left
D-invertibility considers only the membership of a single state to a
single set. For this reason left D-invertibility is much easier to
detect. Our main result (Theorem \ref{Theor D-inv exp-contr}) is a
method to compute left D-invertibility for all linear systems whose
dynamic matrix has no eigenvalue of modulus one.

The main tools used in the paper are the theory of Iterated Function
Systems (IFS), and a theorem of Kronecker. The use of IFS in
relation with left invertibility is described in \cite{Nev}. The
Kronecker's theorem has to do with density in the unit cube of the
fractional part of real numbers. By means of a particular
construction illustrated in section $5$ the problem of ``turning''
left D-invertibility into left invertibility can be handled with a
Kronecker-type density theorem.

The paper is organized as follows. Section $2$ contains an
illustrative example. Section $3$ is devoted to the definitions of
left invertibility and uniform left invertibility. Section $4$
illustrates the background knowledge. In section $5$ left
D-invertibility is introduced and main results (Theorems \ref{Theor
alg independent} and \ref{Theor D-inv exp-contr}) are proved.
Section $6$ contains a deeper study of unidimensional systems. In
section $7$ we present some examples and section $8$ shows
conclusions and future perspectives. The Appendix is devoted to
technical proofs. Finally, in section $9$, we collect the notations
used in the paper.

\section{An illustrative example}

In this section we give an illustrative example, to clarify methods
and the purpose of this paper. Consider the linear output-quantized
system, in one dimension, given by

\begin{equation}\label{eq_syst_illus}
\left\{\begin{array}{l}
  x(k+1) = ax(k)+u(k) \\
  y(k) = \lfloor x(k)\rfloor\\
  u(k)\in \{-1,0,1\},
\end{array}\right.
\end{equation}$\ $\\
where $u$ is the input, $y$ is the output, $a$ is a constant and
$\lfloor\cdot\rfloor$ is the floor function. The left invertibility
problem consists in reconstructing the unknown input sequence $u(k)$
of the system from the information available reading only the
quantized output sequence $y(k)$. For the purpose of this example,
define left invertibility in a negative way (see also definition
\ref{Def_ULI}): if there exist two input sequences
$\{u'(k)\},\{u''(k)\}$ and two initial conditions $x'(0),x''(0)$
such that the resulting orbits $\{x'(k)\}$ and $\{x''(k)\}$ give
rise to the same output, then the system is not invertible (i.e.
there is an output sequence that does not allow the reconstruction
of the input sequence). No matter here how long are the sequences we
are considering. Observe now that $x'(k),x''(k)$ give rise to the
same output if and only if the couple $[x'(k),x''(k)]\in\Bbb R^2$
belongs to the set (see figure 1)

\[Q = \bigcup_{i\in\Bbb Z}\ [i,i+1[\times [i,i+1[\subset\Bbb R^2.\]$\ $\\
Therefore, if there exist two orbits
$\{x'(k)\}_{k=1}^K,\{x''(k)\}_{k=1}^K$ of the system
(\ref{eq_syst_illus}) such that the couple $[x'(k),x''(k)]\in\Bbb
R^2$ remains in $Q$, then the system (\ref{eq_syst_illus}) is not
left invertible. For many reasons (mainly due to the ``complexity''
of the shape of the set $Q$) it is much easier to check the
existence of a sequence of pair of states in $Q_0$, the ``strip''
that includes $Q$ (see figure 1):

\[Q_0 = \left\{ [x',x'']\in\Bbb R^2:\ |x'-x''|<1 \right\}\subset\Bbb
R^2\]$\ $\\
The problem is then to deduce the presence of an orbit in $Q$ from
the existence of an orbit in $Q_0$. We solve this problem with the
help of the following remarks:
\begin{itemize}
  \item[$\spadesuit$] Consider the orbits
  $\{\overline{x(k)'}\}_{k=1}^K$, $\{\overline{x(k)''}\}_{k=1}^K$ of the
  system (\ref{eq_syst_illus}), obtained respectively with initial conditions
  $[x(0)'+\overline{x},x(0)''+\overline{x}]$
  instead of $[x(0)',x(0)'']$. Then easy calculations shows that each $[\overline{x(k)'},\overline{x(k)''}]$
  differs from $[x(k)',x(k)'']$ by a translation along a
  parallel of the bisecting line of $\Bbb R^2$.
  \item[$\clubsuit$] The length of the translation of the state $[\overline{x(k)'},\overline{x(k)''}]$, up to a constant, is
  proportional to $a^k$, since
  \[\overline{x(k)'} = a^k(x(0)'+\overline{x})+\hbox{terms that depend only on the sequence of u's.}\]
  The same holds for $\overline{x(k)''}$.
\end{itemize}
Therefore the question is the following. Consider an orbit
$\Big\{[x(k)',x(k)'']\Big\}$ which is included in $Q_0$. Does there
exist a suitable translation $[\overline{x},\overline{x}]$ of the
initial states $[x(0)',x(0)'']$ such that the resulting orbit is
indeed in $Q$?

By $\spadesuit$, as $\overline x$ varies in $\Bbb R$, the initial
condition is moving along a parallel of the bisector of $\Bbb R^2$,
but every time it cover a distance of $\sqrt 2$, it is in the same
relative position with respect to a square of $Q$: in other words
the property of ``being inside $Q$'' is periodic.

Moreover, by $\clubsuit$, the distance covered by
$[\overline{x(k)'},\overline{x(k)''}]$ with respect to
$[x(k)',x(k)'']$ is proportional to $a^k$, and easy calculations
shows that the property of ``being inside $Q$'', which is periodic,
depends on the fractional parts of $a^k$ (see also figure 2).

\begin{Theorem}
If $1,a,a^2,\ldots,a^K$ are linearly independent over $\Bbb Z$,
then, for every $\alpha_0,\alpha_1,\ldots,\alpha_K\in\Bbb R$ the set
of points
\[\left\{ \left[frac(\alpha_0+l),frac(\alpha_1+la),\ldots,frac(\alpha_K + la^K)\right]: l\in\Bbb R \right\}\]
is dense in the unit cube of $\Bbb R^{K+1}$. $\diamondsuit$
\end{Theorem}

The Theorem (see also Theorem \ref{Theorem Kro 2}) indeed assures
that there exists a suitable translation
$[\overline{x},\overline{x}]$ of the initial states $[x(0)',x(0)'']$
such that the resulting orbit is inside $Q$. Referring to figure 3,
each arrow represents the ``period'' $1$ (in terms of fractional
parts, to be multiplied by $\sqrt 2$), and in particular, the values
of fractional parts going from $0$ to the intersection of the square
with the arrow correspond to a point inside $Q$. Therefore arbitrary
small fractional parts, whose existence is assured by the Theorem,
means that every point of the orbit can be positioned inside $Q$
simply modifying the initial conditions in the way we showed. Indeed
it is well known also that the numbers $a\in\Bbb R$ such that
$1,a,a^2,\ldots,a^K$ are linearly independent over $\Bbb Z$ are a
set of full measure.

That's the point of this paper. In this work we solve the problem of
checking the left invertibility of a system in arbitrary dimension,
by investigating the presence of an orbit in the ``strip'' $Q_0$,
and then deducing the presence of an orbit in $Q$. This strategy
works for a full measure set of parameters.

\begin{figure}\label{fig_ill1}
\begin{center}
\includegraphics[height=3.5in]{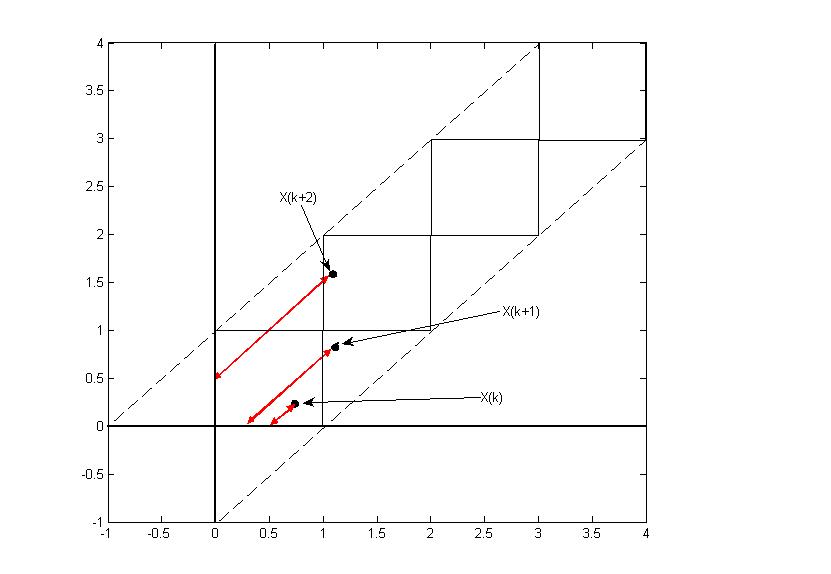}
\caption{The set $Q$ is formed by the squares on the diagonal, while
the set $Q_0$ is the strip inside the dashed line.}
\end{center}
\end{figure}

\begin{figure}\label{fig_ill2}
\begin{center}
\includegraphics[height=3.5in]{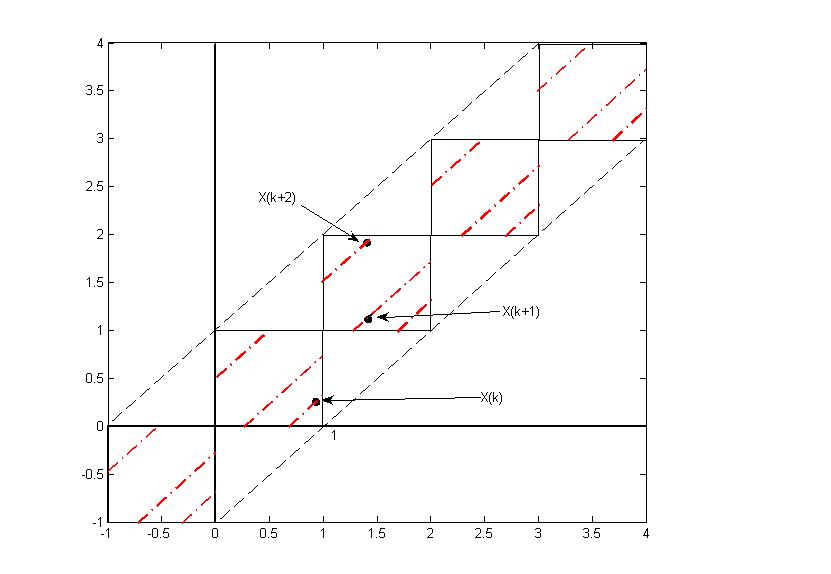}
\caption{The distance covered by
$\Big\{[\overline{x(k)'},\overline{x(k)''}]\Big\}$ with respect to
$\Big\{[x(k)',x(k)'']\Big\}$ is proportional to $a^k$, along a line
parallel to the bisecting line of $\Bbb R^2$. The property of
``being inside $Q$'' is related to fractional parts of $a^k$ which
are associated with dashed line inside the squares.}
\end{center}
\end{figure}

\begin{figure}\label{fig_ill3}
\begin{center}
\includegraphics[height=3.5in]{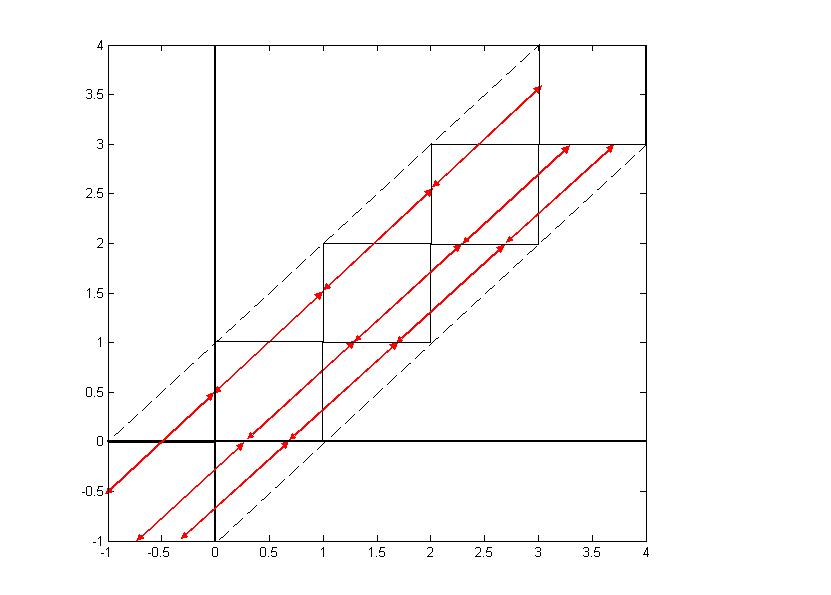}
\caption{Every arrow's length is $\sqrt 2$. The first part of the
arrow, corresponding to small fractional parts, is the part inside a
square of $Q$.}
\end{center}
\end{figure}

\section{Basic setting}

Throughout this paper we use the following notations:
\begin{itemize}
\item $\pi_p$ is the canonical projection on the first $p$
coordinates;

\item $\pi_W$ is the canonical projection on a subspace $W$;

\item $\varpi_i:\Bbb R^d\rightarrow\Bbb R$ denotes the orthogonal
projection on the $i$-th coordinate;

\item $e_i$ is the $i-$th vector of the canonical basis of $\Bbb R^d$;

\item $\lfloor \cdot \rfloor$ is the floor function, acting
componentwise;

\item $frac(\cdot)$ denotes the fractional part, acting
componentwise;

\item $\langle v_1,\ldots,v_i\rangle$ denotes the linear subspace
generated by the vectors $v_1,\ldots,v_i$;

\item $\mu(S)$ denotes the Lebesgue measure of a set $S$ in $\Bbb
R^d$;

\item $\partial$ indicates ``topological boundary of \ldots'';

\item $\setminus$ denotes the set difference;

\item $\Bbb Q[\zeta_1,\ldots,\zeta_N]$ is the ring of polynomial in
the variables $\zeta_1,\ldots,\zeta_N$, with coefficients in $\Bbb
Q$;
\end{itemize}

\begin{Definition}\label{Def_unif_partition}
The uniform partition of rate $\delta$ of $\Bbb R^p$ is
\[\mathcal{P}=\bigcup_{i\in\Bbb Z^p}\mathcal{P}_{i}=
 \bigcup_{i_1,\ldots,i_p\in\Bbb Z}\ [i_1\delta,(i_1+1)\delta[\ \times\ldots\times\ [i_p\delta,(i_p+1)\delta[ ,\]
where $i=i_1,\ldots,i_p.$ $\diamondsuit$
\end{Definition}

We consider systems of the form
\begin{equation}\label{sistema iniziale}
             \left\{\begin{array}{l}
             x(k+1)=Ax(k)+Bu(k)\\
             y(k)=q\big(Cx(k)\big)\\
             \end{array}\right.
\end{equation}
where $A\in\Bbb R^{d\times d},B\in\Bbb R^{m\times d},C\in\Bbb
R^{p\times d}$, $x(k)\in\Bbb{R}^d$ is the state, $y(k)\in\Bbb Z^p$
is the output, and $u(k)\in\mathcal{U}\subset\Bbb R^m$ is the input.
The map $q: \Bbb{R}^p \rightarrow \Bbb Z^p$ is induced by the
uniform partition $\mathcal{P} = \bigcup_{i\in\Bbb Z^p}
\mathcal{P}_i$ of $\Bbb{R}^p$ of rate $\delta$ through $q:(x\in
\mathcal P_i) \mapsto i$ and will be referred to as the output
quantizer. We assume that $\mathcal{U}$ is a finite set of
cardinality $n$.

\begin{Remark}\label{rem_canon_syst}
Suitably changing bases, without loss of generality in the system
(\ref{sistema iniziale}) we can suppose $\delta=1$ and $C=\pi_p$.
$\diamondsuit$
\end{Remark}

So we consider only systems of the form
\begin{equation}\label{sistema}
    \left\{\begin{array}{l}
             x(k+1)=Ax(k)+Bu(k)=f_{u(k)}(x(k))\\
             y(k)=\lfloor\pi_px(k)\rfloor.
           \end{array}
    \right.
\end{equation}

\begin{Definition}
A pair of input strings $\{u(i)\}_{i\in\Bbb{N}}$,
$\{u^{\prime}(i)\}_{i\in\Bbb{N}}$ is uniformly distinguishable in
$k$ steps if there exists $l\in\Bbb N$ such that $\forall x(0),
x^{\prime}(0) \in \Bbb R^d$ and $\forall m>l$ the following holds
for the correspondent orbits:
\[u(m)\not=u^{\prime}(m)\;\Rightarrow\;[y(m+1),\ldots,y(m+k)] \not=
[y^{\prime}(m+1),\ldots,y^{\prime}(m+k)],\] (outputs $y(i)$ are
referred to the system with initial condition $x(0)$ and inputs
$u(i)$, while outputs $y^{\prime}(i)$ are referred to the system
with initial condition $x^{\prime}(0)$ and inputs $u^{\prime}(i)$).
In this case, we say that the strings are uniformly distinguishable
with waiting time $l$. $\diamondsuit$
\end{Definition}

\begin{Definition}\label{Def_ULI}
A system of type (\ref{sistema}) is uniformly left invertible (ULI)
in $k$ steps if every pair of distinct input sequences is uniformly
distinguishable in $k$ steps after a finite time $l$, where $k$ and
$l$ are constant. $\diamondsuit$
\end{Definition}

For a ULI system, it is possible to recover the input string until
instant $m$ observing the output string until instant $m+k$. For
applications, it is important to obtain an algorithm to reconstruct
the input symbol used at time $m>l$ by processing the output symbols
from time $m$ to $m+k$.\\

\begin{Definition}\label{Def set Q}
Define the quantization set relative to the system (\ref{sistema})
to be
\[Q = \bigcup_{i_1,\ldots i_p\in\Bbb Z}\underbrace{\Big\{ [i_1,i_{1}+1[\times\ldots
[i_p,i_{p}+1[\times \langle e_{p+1},\ldots,e_d\rangle
\Big\}}_{\subset \langle e_1,\ldots,e_d\rangle}\ \times\ \]
\[ \times\underbrace{\Big\{ [i_1,i_{1}+1[\times\ldots [i_p,i_{p}+1[\times
\langle e_{d+p+1},\ldots,e_{2d}\rangle \Big\}}_{\subset \langle
e_{d+1},\ldots,e_{2d}\rangle}\subset\Bbb{R}^{2d}\] i.e. $Q$ contains
all pairs of states that are in the same element of the partition
$\mathcal{P}$. $\diamondsuit$
\end{Definition}

To address left invertibility, we are interested in the following
system on $\Bbb{R}^{2d}$:

\begin{Definition}
Define the doubled system relative to the system (\ref{sistema}) to
be

\begin{equation}\label{Sistema2d}
X(k+1) = \left[\begin{array}{c}
                    Ax_1(k)+Bu(k) \\
                    Ax_2(k)+Bu^{\prime}(k)
                \end{array}\right]
\end{equation}$\ $\\
where $X(k) = \left(\begin{array}{c} x_1(k) \\
                                     x_2(k)
                     \end{array}\right), \ \ U(k) = \left(\begin{array}{c} u(k) \\
                                                                          u^{\prime}(k)
                                                         \end{array}\right).$
$\diamondsuit$
\end{Definition}

If there exist sequences $\{u(k)\},\{u^{\prime}(k)\}$, and an
initial state in $Q$ such that the corresponding orbit of
(\ref{Sistema2d}) remains in $Q$, then the two strings of inputs
generate the same output for the system (\ref{sistema}). So
conditions ensuring that the state is outside $Q$ for some $k$ will
be investigated to guarantee left invertibility.

\section{Background: attractors and left invertibility}

In this section we recall some results coming from Iterated Function
System theory (see \cite{Bar,FG} for general theory about IFS), in
connection with the notions of left invertibility.

\begin{Definition}
An output-quantized linear system of type (\ref{sistema}) is joint
contractive if $|\lambda|<1$ for every eigenvalue $\lambda$ of the
matrix $A$. It is joint expansive if $|\lambda|>1$ for every
eigenvalue $\lambda$ of the matrix $A$. $\diamondsuit$
\end{Definition}

\begin{Definition}\label{Def_attractor}
Consider an output-quantized system.
\begin{itemize}
  \item A set $\mathcal{A}$ is an attractor if for all orbits $\{x(k)\}_{k\in\Bbb
  N}$ it holds

  \[\lim_{k\rightarrow\infty} dist\big(x(k),\mathcal{A}\big)=0.\]$\ $\\
  Here $dist(x,\mathcal A)$ is the inf of distances between $x$ and points of $\mathcal A$.
  \item A set $\mathcal{I}$ is an invariant set if

  \[\mathcal{I} = \bigcup_{u\in\mathcal U} A(\mathcal{I})+Bu.\ \diamondsuit\]
\end{itemize}
\end{Definition}

\begin{Theorem}\label{Teorema insieme fisso}\cite{Bar,Gw}
Let a system be joint contractive. Then, for every $u\in\mathcal
U^\Bbb N$ the limit $\phi(u) =
\lim_{k\rightarrow\infty}A^kx+A^{k-1}u(k-1)+\ldots+Bu(0)$ exists for
every $x$ and is independent of $x$. The set $\phi(\mathcal U^\Bbb
N)$ is the unique compact attractor and invariant set of the system.
$\diamondsuit$
\end{Theorem}$\ $\\

Consider now an output-quantized linear system, together with two
sets, namely $\mathcal H$ and $S$: $\mathcal H$ plays the role of an
attractive and invariant set, in which the dynamic of the system is
confined, and $S$ plays the role of a quantization set, i.e. a set
that the state has to exit to guarantee an invertibility property.
In \cite{Nev2} a necessary and sufficient condition for left
invertibility of joint contractive systems is given, but here we
state the same condition in a more abstract setting: $\mathcal H$
and $S$ are an attractor and a quantization set, not \emph{the}
attractor and \emph{the} quantization set of the system
(\ref{sistema}). That's because in the following we will use these
results for another attractor and quantization set, i.e. those ones
of the difference systems.

\begin{Definition}\label{Def. Graph associaed to attractor}
The graph $G_k$ associated to the attractor $\mathcal{H}$ is given by:\\
$\bullet$ The set of vertices
\[V = \{\mathcal{H}_{u(0)\ldots u(k)} = A^{k+1}(\mathcal
H)+A^{k}Bu(0)+\ldots+Bu(k):
  u(i)\in\mathcal U\}.\]
$\bullet$ There is an edge from
  $\mathcal{H}_{u(0)\ldots u(k)}$ to
  $\mathcal{H}_{u^{\prime}(0)\ldots u^{\prime}(k)}$ if and only if
  $u(i+1)=u^{\prime}(i)$, for $i=0,\ldots,k-1$. In this case we say
  that the edge is induced by the input $u^{\prime}(k)$. $\diamondsuit$
\end{Definition}

\begin{Definition}\label{Def Int Ext Inv Graph}
Consider the graph $G_k$, and delete all vertices (together with all
starting and arriving edges) $\mathcal{H}_{u(0)\ldots u(k)}$ such
that $\mathcal{H}_{u(0)\ldots u(k)}\cap \{\Bbb{R}^{d}\setminus
S\}\not=\emptyset$. This new graph is called internal invertibility
graph, and denoted with $IG_k$. The union of vertices (which are
sets) of $IG_k$ is denoted by $V_{IG_k}.$
\end{Definition}

\begin{Theorem}\cite{Nev}\label{Theorem Invert iff}
Denote with $\partial S$ the boundary of $S$. Suppose that
$\mathcal{H}\cap \partial S=\emptyset$. Then there exists a
(computable) $k$ such that $V_{IG_k} = V_{IG_k}\cap S$.
$\diamondsuit$
\end{Theorem}$\ $\\

If instead the system (\ref{sistema}) is joint expansive, then the
map $x(k)\mapsto x(k+1)$ admits an inverse for every
$u\in\mathcal{U}$. Therefore it is possible to define a
correspondent inverse system:

\begin{Definition}
If the system (\ref{sistema}) is joint expansive, the inverse system
is

\[\tilde x(k+1) = A^{-1}[\tilde x(k)-Bu(k)]\].$\ $\\
The inverse doubled system, relative to the doubled system
(\ref{Sistema2d}), is defined in a similar way. $\diamondsuit$
\end{Definition}

In the case of joint contractive systems, the inverse systems give
rise to attractors, since they are joint contractive: such
attractors can be described also as the set of initial conditions
that can start a bounded orbit of the system (\ref{sistema}) or of
the doubled system (\ref{Sistema2d}):

\begin{Theorem}\label{Th cond nec joint exp}\cite{Nev}
Suppose that the system (\ref{sistema}) is joint expansive. If there
is an infinite bounded orbit of the system (\ref{sistema}) or of
doubled system (\ref{Sistema2d}), then this orbit is entirely
contained in the attractor of the inverse system or in the attractor
of the inverse doubled system, respectively. Consequently, if we
restrict to bounded orbits, Theorem \ref{Theorem Invert iff} applies
to these attractors. $\diamondsuit$
\end{Theorem}

\section{Difference system and D-invertibility}

\begin{Definition}\label{Def Difference System}
The difference system associated with the system (\ref{sistema}) is
\begin{equation}\label{equation difference set}
    z(k+1) = Az(k)+Bv(k)
\end{equation}
where $z(k)\in\Bbb R^d$,
$v(k)\in\mathcal{V}=\mathcal{U}-\mathcal{U}=\{u-u^{\prime}:\
u\in\mathcal U,\ u^{\prime}\in\mathcal U\}$. $\diamondsuit$
\end{Definition}

\begin{Remark}\label{Remark difference system leq 1}
The difference system represents at any instant the difference
between the two states $x(k)-x^{\prime}(k)=z(k)$ when the input
symbols $u(k)-u^{\prime}(k)=v(k)$ are performed. So we are
interested in understanding the conditions under which

\[\{z(k)\}\ \cap\ \{\ ]-1,1[\ \}^p\times \langle e_{p+1},\ldots,e_{d}\rangle =
\emptyset.\]$\ $\\
Indeed, this implies that $y(k)\not=y^{\prime}(k)$. The converse is
obviously not true. $\diamondsuit$
\end{Remark}

\begin{Definition}\label{Def_Dk1k2}
Consider the difference system. If $z(0)$ is an initial condition
and $(v(1),\ldots,v(k_2))$ a sequence of inputs of the difference
system, we let $D^{k_2}_{k_1}(z(0),v(1),\ldots,v(k_2))$ denote the
sequence $\left(\pi_pz(k_1),\ldots,\pi_pz(k_2)\right)$ generated by
the difference system (\ref{equation difference set}) with initial
condition $z(0)$ and input string $(v(1),\ldots,v(k_2))$.
$\diamondsuit$
\end{Definition}

\begin{Definition}\label{Def_unif_disting}
A pair of input strings $\{u(i)\}_{i\in\Bbb{N}}$,
$\{u^{\prime}(i)\}_{i\in\Bbb{N}}$ is uniformly D-distinguishable in
$k$ steps if there exists $l\in\Bbb N$ such that $\forall x(0),
x^{\prime}(0) \in \Bbb R^d$ and $\forall m>l$ the following holds:

\[v(m)\not=0\ \Rightarrow\ D_{m+1}^{m+k}(z(0),v(1),\ldots,v(m+k))\not\in\ \underbrace{]-1,1[^p\ \times\ldots\times\ ]-1,1[^p}_{k\
times},\]$\ $\\
where $z(0)=x(0)-x^{\prime}(0)$ and $v(i)=u(i)-u^{\prime}(i)$. In
this case, we say that the strings are uniformly D-distinguishable
with waiting time $l$. $\diamondsuit$
\end{Definition}

\begin{Definition}\label{def_ULDI}
A system of type (\ref{sistema}) is uniformly left D-invertible
(ULDI) in $k$ steps if every pair of distinct input sequences is
uniformly D-distinguishable in $k$ steps after a finite time $l$,
where $k$ and $l$ are constant. $\diamondsuit$
\end{Definition}

\begin{Remark}
Thanks to Remark \ref{Remark difference system leq 1} uniform left
D-invertibility implies uniform left invertibility. $\diamondsuit$
\end{Remark}

The first main result is based on a density theorem of Kronecker.

\begin{Definition}\label{def_lin_independent}
The numbers $\vartheta_1,\ldots,\vartheta_M\in\Bbb R$ are linearly
independent over $\Bbb Z$ if the following holds:
\[ k_1,\ldots,k_M\in\Bbb Z: k_1\vartheta_1+\ldots,+k_M\vartheta_M=0\ \ \Rightarrow\ \ k_1=\ldots=k_M=0.\ \diamondsuit\]
\end{Definition}

\begin{Theorem}[Kronecker]\cite{Har}\label{Theorem Kro 2}
If $\vartheta_1,\ldots,\vartheta_M,1\in\Bbb R$ are linearly
independent over $\Bbb Z$, then, for every
$\alpha_0,\alpha_1,\ldots,\alpha_M\in\Bbb R$ the set of points
\[\left\{ \left[frac(\alpha_0+l),frac(\alpha_1+l\vartheta_1),\ldots,frac(\alpha_M+l\vartheta_M)\right]: l\in\Bbb R \right\}\]
is dense in the unit cube of $\Bbb R^{M+1}$. $\diamondsuit$
\end{Theorem}

Considering the difference system (Definition \ref{Def Difference
System}), we are interested in orbits completely included in
$(]-1,1[)^p\times\langle e_{p+1},\ldots,e_{d}\rangle.$ The following
proposition shows that under a very weak condition orbits completely
included in $(]-1,1[)^p\times\langle e_{p+1},\ldots,e_{d}\rangle$
must be bounded.

\begin{Proposition}\label{Prop transverse bounded}
Suppose that the matrix $A$ does not have an invariant subspace
included in $\langle e_{p+1},\ldots,e_{d}\rangle$. Then there exists
a bounded set $I$ such that, if $\{z(k)\}_{k\in\Bbb N}\subset
(]-1,1[)^p\times\langle e_{p+1},\ldots,e_{d}\rangle$ is an orbit of
the difference system, then $\{z(k)\}_{k\in\Bbb N}\subset I$.
$\diamondsuit$
\end{Proposition}

\emph{Proof:} See Appendix. $\diamondsuit$$\ $\\

Note that the set of matrices $A\in\Bbb R^{d\times d}$ that have an
invariant subspace in $\langle e_{p+1},\ldots,e_{d}\rangle$ is a
zero measure set. Define now $\Bbb S_D(B,\mathcal{U})$ to be the set
of matrices $A\in\Bbb R^{d\times d}$ such that the system
(\ref{sistema}) is uniformly left D-invertible, and $\Bbb
S(B,\mathcal{U})$ to be the set of matrices $A\in\Bbb R^{d\times d}$
such that the system (\ref{sistema}) uniformly left invertible.

\begin{Definition}\label{def_alg_independent}
Indicate with $\Bbb Q[\zeta_1,\ldots,\zeta_N]$ the ring of
polynomials in the variables $\zeta_i$ with coefficients in $\Bbb
Q$. The set of numbers $\alpha_1,\ldots,\alpha_N\in\Bbb C$ is said
to be algebraically independent if
\[0\not=p(\zeta_1,\ldots,\zeta_N)\in\Bbb Q[\zeta_1,\ldots,\zeta_N]\ \Rightarrow\ p(\alpha_1,\ldots,\alpha_N)\not=0.\ \diamondsuit\]
\end{Definition}

\begin{Theorem}\label{Theor alg independent}
Suppose that in the system (\ref{sistema}) the set of elements of
the matrix $A$ is algebraically independent. Then the system is
uniformly left D-invertible if and only if it is uniformly left
invertible. This in turn implies that $\Bbb
S(B,\mathcal{U})\setminus\Bbb S_D(B,\mathcal{U)}$ has measure zero
in $\Bbb R^{d\times d}$ for every $B,\mathcal{U}$.
\end{Theorem}

\emph{Proof:} See appendix. $\diamondsuit$

\subsection{D-invertibility of output-quantized linear systems}

We are going to show how to detect left D-invertibility of any
linear systems without eigenvalues of modulus one. Suppose that, if
$\lambda$ is an eigenvalue of the matrix $A$, then
$|\lambda|\not=1$. Denote with $E_c,E_e$ respectively the
contractive and the expansive eigenspaces of the matrix $A$, i.e.
the eigenspaces relative to eigenvalues $<1$ and $>1$ in modulus,
respectively. Because of the hypothesis on the eigenvalues we have
$E_c + E_e=\{x+y:\ x\in E_c,\ y\in E_e\} = \Bbb R^d$. Now consider
the following two systems respectively on $E_c,E_e$, that are joint
contractive:

\[z_c(k+1) = (\pi_{E_c}A)z_c(k)+\pi_{E_c}(Bv(k)));\]

\[z_e(k+1) = (\pi_{E_e}A^{-1})[z_e(k)-\pi_{E_e}(Bv(k))],\]$\ $\\
where with $z_c,z_e$ we indicate the projections of $z$ onto
$E_c,E_e$, respectively. $v(k)\in\mathcal V$. The above systems must
have invariant attractors $\mathcal{T}^c,\mathcal{T}^e$. Let us
denote with $\mathcal{T}$ the attractor
\[\mathcal T = \mathcal{T}^c+\mathcal{T}^e = \{x+y:\ x\in \mathcal T^c,y\in\mathcal T^e\}.\]
We can now apply the construction of internal invertibility graph
(Definition \ref{Def Int Ext Inv Graph}) for the attractor
$\mathcal{T}$ (i.e. substituting $\mathcal H$ with $\mathcal T$),
substituting $\mathcal S$ with $\big\{\ ]-1,1[\
\big\}^p\times\langle e_{p+1},\ldots,e_{d}\rangle$, and calling a
path $\{V_1,\ldots,V_i\}$ on $IG_k$ proper if it is induced by an
input $v\in\mathcal{V},\ v\not=0$. Denoting with
$\partial\Big((]-1,1[)^p\times\langle
e_{p+1},\ldots,e_{d}\rangle\Big)$ the boundary of
$(]-1,1[)^p\times\langle e_{p+1},\ldots,e_{d}\rangle$ we have the
following

\begin{Theorem}\label{Theor D-inv exp-contr}
Suppose that $|\lambda|\not=1$ for every eigenvalue $\lambda$ of the
matrix $A$, that $\mathcal{T}\cap
\partial\Big((]-1,1[)^p\times\langle
e_{p+1},\ldots,e_{d}\rangle\Big)=\emptyset$, and that $A$ does not
have an invariant subspace in $\langle e_{p+1},\ldots,e_{d}\rangle$.
Then the system (\ref{sistema}) is uniformly left D-invertible if
and only if $IG_k$ does not contain arbitrary long proper paths,
where $k$ is the one identified through Theorem \ref{Theorem Invert
iff}.
\end{Theorem}

\emph{Proof:} Since $A$ does not have an invariant subspace in
$\langle e_{p+1},\ldots,e_d\rangle$, by Proposition \ref{Prop
transverse bounded} all orbits of the difference system included in
$(]-1,1[)^p\times\langle e_{p+1},\ldots,e_{d}\rangle$ must be
bounded, and, by Theorem \ref{Th cond nec joint exp}, must be
included in the attractor $\mathcal{T}$. By Theorem \ref{Theorem
Invert iff} the system is uniformly D-invertible if and only if
$IG_k$ does not contain arbitrary long proper paths. $\diamondsuit$

\begin{Remark}
Theorem \ref{Theor D-inv exp-contr} gives an explicit,
algorithmically implementable, way to compute ULDI of a system. By
Theorem \ref{Theor alg independent} we are able to compute ULI in
the same way for systems with a full measure set of matrices $A$
(the conditions on the eigenspaces do not affect the full measure).
$\diamondsuit$
\end{Remark}

\begin{Remark}
The technical condition $\mathcal{T}\cap
\partial\Big((]-1,1[)^p\times\langle
e_{p+1},\ldots,e_{d}\rangle\Big)=\emptyset$ means, from a practical
point of view, that left D-invertibility can be checked up to any
finite precision, since the set $\Big((]-1,1[)^p\times\langle
e_{p+1},\ldots,e_{d}\rangle\Big)$ depends only on the partition
$\mathcal P$, and any small ``disturbance'' of the rate of the
partition $\mathcal P$ allows the application of the Theorem
\ref{Theor D-inv exp-contr}. Further details on this point are given
in \cite{Nev2}. $\diamondsuit$
\end{Remark}

\section{Output-quantized linear systems of dimension $1$}

Linear systems of dimension $1$ assume the following form, deriving
from (\ref{sistema}):

\begin{equation}\label{Eq unid lin systems}\left\{\begin{array}{l}
x(k+1) = ax(k)+u(k)\\
y(k) = \lfloor x(k)\rfloor\\
\end{array}\right.\end{equation}
This is a contractive system if $|a|<1$ and an expansive system if
$|a|>1$. If $|a|<1$ the invertibility problem can be solved with the
methods of section $3$ (see \cite{Nev}). The next Theorem shows a
necessary condition for the ULI of a system of type (\ref{Eq unid
lin systems}): if it is not satisfied we construct inductively a
pair of strings that gives rise to the same output.

\begin{Theorem}\label{Controesempio}
Suppose that in the system (\ref{Eq unid lin systems}) $|a|>2$. If
there exist $u_1,u_2\in\mathcal{U},u_1\not=u_2$ such that
$|u_1-u_2|<a$, then the system is not ULI.
\end{Theorem}

\emph{Proof:} We will consider sequences of sets of type
\begin{equation}\label{eq set succession}\left\{\begin{array}{l}
           S_{i+1}= \left\{a(S_{i})+u(i)\right\}\bigcap\left\{a(S_{i})+u^{\prime}(i)\right\}\bigcap \mathcal{P}(i+1) \\
           S_0 = [0,1[,
         \end{array}\right.\end{equation}
where $u(i),u^{\prime}(i)\in\{u_1,u_2\}$ and
$\mathcal{P}(i+1)\in\mathcal{P}$ is chosen at each step to maximize
the measure of $S_{i+1}$.

In the sequence (\ref{eq set succession}) set $u(1)=u_1$,
$u^{\prime}(1)=u_2$. Since $|u_1-u_2|<a$, there exists a
$\mathcal{P}(1)\in\mathcal{P}$ such that $\mu(S_1)>0$ (recall that
$\mu$ indicates the Lebesgue measure). Then, for $i>1$ define
\[u(i)=u^{\prime}(i)=u_1.\]
Since $|a|>2$ there exists an $i_0$ such that
$\mu\left(S_{i_0}\right)=1$, therefore, applying again
$u(i_0+1)=u_1$ and $u^{\prime}(i_0+1)=u_2$
\[\mu\left\{A(S_{i_0})+Bu_1\ \cap\ A(S_{i_0})+Bu_2\right\}>0.\]
So there exists $x_0,x^{\prime}_0\in\Bbb R$ and
$(u(1),\ldots,u(i_0+1)),(u^{\prime}(1),\ldots,u^{\prime}(i_0+1)),$
with $u(1)\not=u^{\prime}(1)$ and $u(i_0+1)\not=u^{\prime}(i_0+1)$,
such that for the corresponding outputs it holds
\[(y(0),\ldots,y(i_0+1)) = (y^{\prime}(0),\ldots,y^{\prime}(i_0+1))\]

It is then enough to point out that, since we can achieve every pair
of states $x,x'\in S_{i_0}$ in the above described way, we can again
go on in the same way and find a new instant $i_1$, a pair of
initial states $x_{1,0},x'_{1,0}$, and control sequences
$(u(1),\ldots,u(i_1)),\ (u^{\prime}(1),\ldots,u'(i_1))$, with
$u(i_1)\not=u^{\prime}(i_1)$, such that for the corresponding output
it holds
\[(y(0),\ldots,y(i_1)) = (y^{\prime}(0),\ldots,y^{\prime}(i_1)).\]

Finally, we can achieve by induction an increasing finite sequence,
but arbitrarily long, of instants $i_k$, pairs of initial states
$(x_{k,0},x'_{k,0})$, and sequences of controls
$(u(1),...,u(i_k)),(u'(1),...,u'(i_k))$ with $u(i)\not =u'(i)$ if
$i=i_j+1$ for $j=1,...,k-1$ such that such that for the
corresponding output it holds
\[(y(0),\ldots,y(i_k)) = (y^{\prime}(0),\ldots,y^{\prime}(i_k)).\]
This contradicts the uniform left invertibility property.
$\diamondsuit$

\begin{Definition}\label{def_alg_trascendental}
A number $\alpha\in\Bbb C$ is called algebraic if there exists a
polynomial $p(x)\in\Bbb Z[x]$ such that $p(\alpha)=0$. In this case
the minimum degree of a polynomial with such a property is called
the degree of $\alpha$. A number $\alpha\in\Bbb C$ is called
trascendental if it is not algebraic. $\diamondsuit$
\end{Definition}

The following Theorem can be deduced from Theorem \ref{Theor alg
independent}, observing that an algebraically independent set of one
element is a trascendental number.

\begin{Theorem}
Suppose that $a$ is trascendental. Then the system (\ref{Eq unid lin
systems}) is uniformly invertible if and only if it is uniformly
D-invertible. $\diamondsuit$
\end{Theorem}

\begin{Proposition}
The unidimensional system (\ref{Eq unid lin systems}) is either ULDI
in time 1, or not ULDI at all.
\end{Proposition}

\emph{Proof:} A sufficient condition for uniform left
D-invertibility in one step is
\[\forall v\in\mathcal{V},v\not=0:\ |v|\geq |a|+1:\]
indeed in this hypothesis $\forall v\in\mathcal{V},v\not=0$
\[ ]-1,1[\ \cap\ \big\{a\cdot(]-1,1[)+v\big\} =\ ]-1,1[\ \cap\ ]-a+v,a+v[\ = \emptyset\]

We now prove that if $\exists v\in\mathcal{V},v\not=0:\ |v|< |a|+1$,
then the system is not uniformly left D-invertible. Indeed in this
case the system
\[\left\{\begin{array}{c}
    ax_1+v=x_2 \\
    ax_2-v=x_1
  \end{array}\right.\]
has the solution $x_1=\frac{-v}{a+1},x_2=\frac{v}{a+1}$. Since
$|x_1|,|x_2|<1$ the difference system has the infinite orbit
$\{x_1,x_2,x_1,x_2,\ldots\}\subset (]-1,1[)^p\times\langle
e_{p+1},\ldots,e_{d}\rangle$. Therefore system (\ref{Eq unid lin
systems}) is not left D-invertible. $\diamondsuit$

\begin{Corollary}
Consider the unidimensional system (\ref{Eq unid lin systems}), with
trascendental $a$. Then it is either ULI in one step, or it is not
ULI. $\diamondsuit$
\end{Corollary}

\begin{Remark}
It's easy to see that a system of the form (\ref{Eq unid lin
systems}) is uniformly D-invertible in one step if for all
$u_1,u_2\in\mathcal{U}$ it holds $|u_1-u_2|>|a|+1$. Therefore we
have this summarizing situation for unidimensional systems:
\begin{itemize}
  \item[$\bullet$] $|a|<1$: ULI can be detected with methods described in
  section $3$.
  \item[$\bullet$] $|a|>2$: $\left\{
                   \begin{array}{ll}
                     \min_{u_1,u_2\in\mathcal{U}}|u_1-u_2|<|a|, & \hbox{the system is not ULI;} \\
                     \min_{u_1,u_2\in\mathcal{U}}|u_1-u_2|>|a|+1, & \hbox{the system is ULI in 1 step;} \\
                     |a|\leq\min_{u_1,u_2\in\mathcal{U}}|u_1-u_2|\leq |a|+1, & \hbox{$\bigstar$ holds.}
                   \end{array}
                 \right.$
  \item[$\bullet$] $1\leq |a|\leq 2$: $\bigstar$ holds.
  \item[$\bigstar$] The unidimensional system (\ref{Eq unid lin systems}) is either ULDI
in time 1, or not ULDI at all. With the additional hypothesis of
trascendence of $a$ the system is ULI in one step or it is not ULI.
$\diamondsuit$
\end{itemize}
\end{Remark}

\section{Examples}

Consider the system (\ref{sistema}). As stated in Theorem \ref{Theor
alg independent}, if the elements of the matrix $A$ forms an
algebraically independent set, then uniform left D-invertibility is
equivalent to uniform left invertibility. A standard method to
construct algebraically independent sets can be easily deduced from
the following Theorem of Lindemann and Weierstrass:

\begin{Theorem}\cite{Bak}\label{Theor Lin-Weier}
Suppose that the numbers $\alpha_1,\ldots,\alpha_N$ are linearly
independent over $\Bbb Q$. Then $e^{\alpha_1},\ldots,e^{\alpha_N}$
are an algebraically independent set. $\diamondsuit$
\end{Theorem}

\begin{Example} Consider the system (\ref{sistema}) with
\begin{equation}\label{eq example 1 system}
A=\left(
    \begin{array}{cc}
      e^{\sqrt{5}} & e^{\sqrt{3}} \\
      e^{\sqrt{2}} & e^{\sqrt{7}} \\
    \end{array}
  \right),\ \ B = \left(
                 \begin{array}{c}
                   1 \\
                   1 \\
                 \end{array}
               \right),\ \ \mathcal{U}=\{0,\pm 1, 2\}\ \ y(k)=\lfloor\pi_1x(k)\rfloor.
\end{equation}
The two eigenvalue of $A$ are approximately $6.3531$ and $17.0974$,
so system (\ref{eq example 1 system}) is joint expansive. The
difference system is then joint expansive too. The inverse
difference system is given by
\begin{equation}\label{eq inv diff syst ex1}
x(k+1)=A^{-1}(x(k)-v(k)),\ \ v(k)\in\mathcal{V}=\{0,\pm 1,\pm 2,\pm
3\}
\end{equation}
which is joint contractive. It is possible to show (see figure 4)
that the attractor of the inverse difference system (\ref{eq inv
diff syst ex1}) is included in $]-1,1[\times\Bbb R.$ So the system
is not uniformly left D-invertible.

Moreover the elements of the matrix $A$, by Theorem \ref{Theor
Lin-Weier}, are an algebraically independent set because
$\sqrt{5},\sqrt{2},\sqrt{3},\sqrt{7}$ are linearly independent over
$\Bbb Q$. So by Theorem \ref{Theor alg independent} system (\ref{eq
example 1 system}) is not uniformly left-invertible. $\diamondsuit$
\end{Example}

\begin{figure}\label{FigExample1}
\begin{center}
\includegraphics[width=0.7\textwidth]{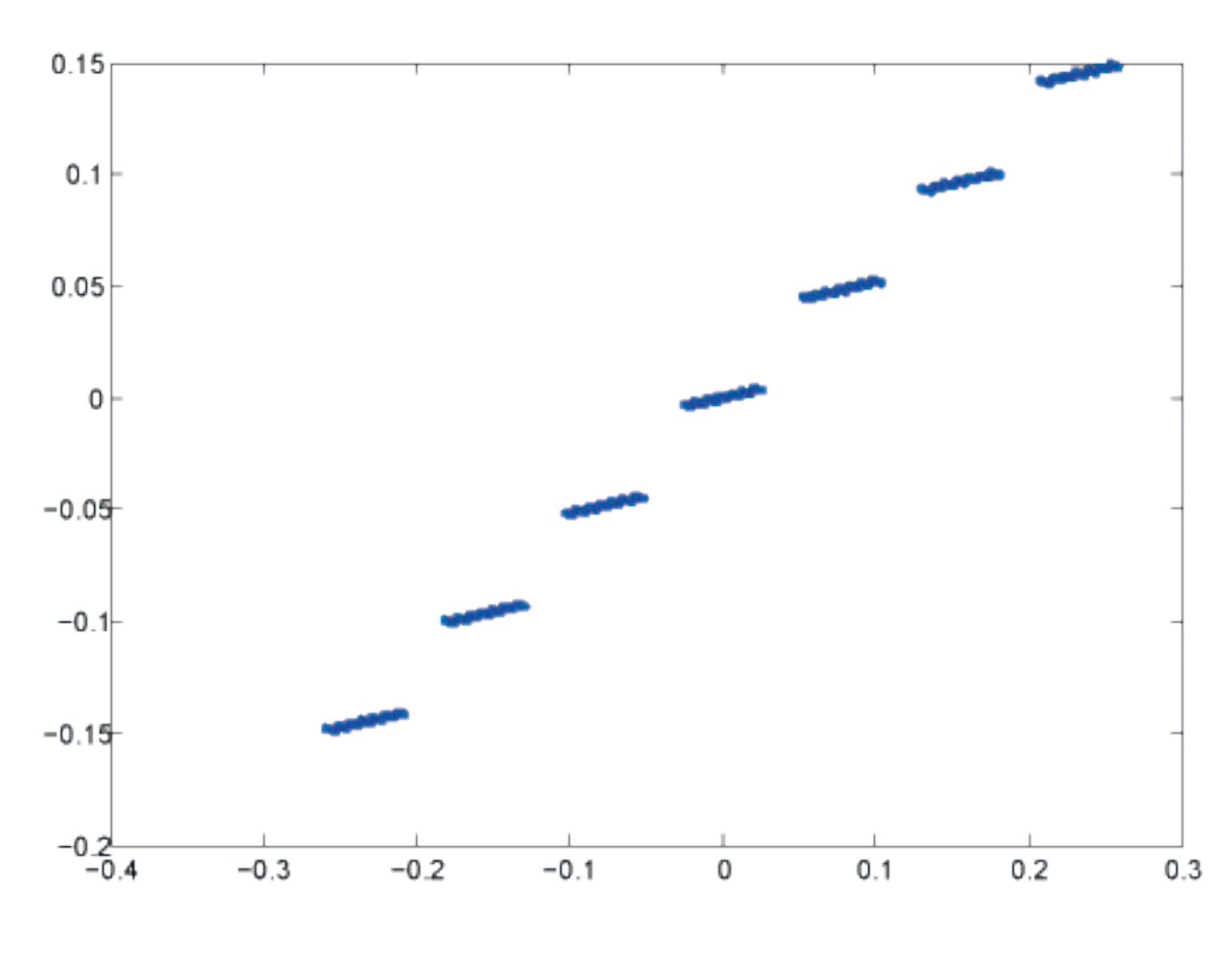}
\caption{Attractor of system (\ref{eq inv diff syst ex1}).}
\end{center}
\end{figure}

\begin{Example}
Consider the system (\ref{sistema}) with
\begin{equation}\label{eq example 2 system}
A=\left(
    \begin{array}{ccc}
      2 &    0        & 0 \\
      0 & \frac{1}{2} & 1 \\
      0 &    0        & \frac{1}{3} \\
    \end{array}
  \right),\ \ B = \left(
                 \begin{array}{c}
                   2 \\
                   3 \\
                   6 \\
                 \end{array}
               \right),\ \ \mathcal{U}=\{0, 1\}\ \ y(k)=\lfloor\pi_2x(k)\rfloor.
\end{equation}
We have so $\mathcal{V}=\{0,\pm 1\}$. The three eigenvalue of the
matrix are $\frac{1}{2}$, $\frac{1}{3}$ and $2$, so we can apply
Theorem \ref{Theor D-inv exp-contr}. Therefore we split $\Bbb R^3$
in $E_c=\Bbb R^2$ (identified with $\{0\}\times\Bbb R^2$) and
$E_e=\Bbb R$ (identified with $\Bbb R\times\{0\}\times\{0\}$). The
attractor $\mathcal{A}_e$ relative to the inverse difference system
on $E_e$ is $[-2,2]$, while the attractor $\mathcal{A}_c$ relative
to the inverse difference system on $E_c$ is drawn in figure 5. We
are interested in orbits of the inverse difference system on $E_e$
that remains in $]-1,1[$, and in orbits of the difference system on
$E_c$ that remains in $]-1,1[\times\Bbb R$. It's easy to see that

\[\left(
    \begin{array}{cc}
      \frac{1}{2} & 1 \\
      0 & \frac{1}{3} \\
    \end{array}
    \right)\left(\mathcal{A}_c\cap\ ]-1,1[\times\Bbb R\right)+\left(
                                                    \begin{array}{c}
                                                    3 \\
                                                    6 \\
                                                    \end{array}
                                                  \right)
\ \bigcap\ \big\{\mathcal{A}_c\cap\ ]-1,1[\ \times\Bbb R\big\} =
\emptyset,\]$\ $\\
so, no matter the behavior of the system on $E_e$, system (\ref{eq
example 2 system}) is uniformly left D-invertible in one step,
therefore uniformly left invertible in one step. $\diamondsuit$
\end{Example}

\begin{figure}\label{FigExample2}
\begin{center}
\includegraphics[width=0.7\textwidth]{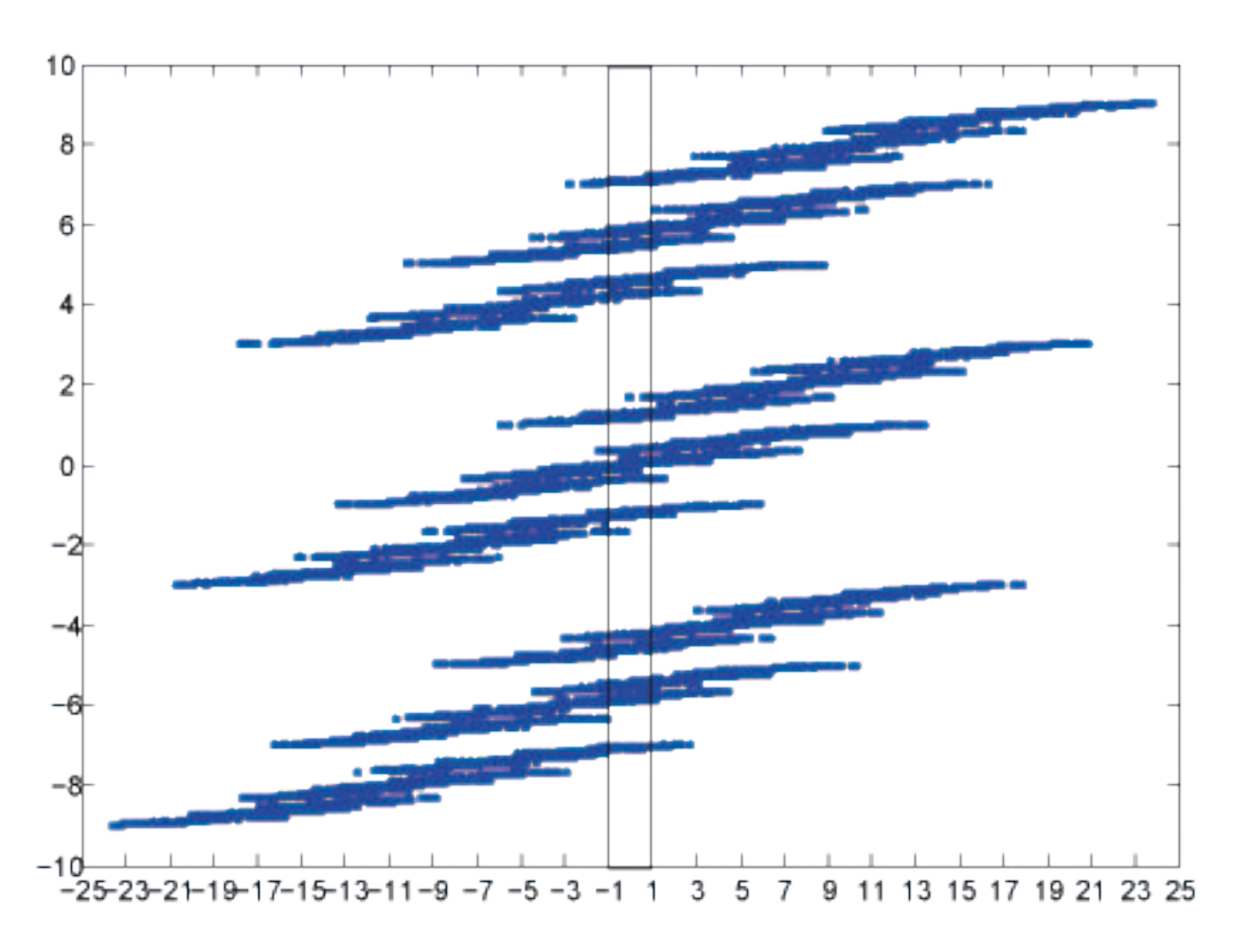}
\caption{Attractor of the inverse difference system on $E_c$.}
\end{center}
\end{figure}

The last example illustrates the difference between left
D-invertibility and left invertibility.

\begin{Example}
Consider the unidimensional system (\ref{Eq unid lin systems}) with
\begin{equation}\label{eq example 3 system}
a=1/2,\ \ \ \mathcal U=\{-1,0,1\}.
\end{equation}
We are going to show that system (\ref{eq example 3 system}) is
uniformly left invertible but not uniformly left D-invertible.

To show that the system is not ULDI consider the following orbit
with initial condition $x_0=\frac{1}{2}$:

\[x_{k+1} = \left\{
              \begin{array}{ll}
                \frac{1}{2}x(k)+1, & \ \ \ \ \hbox{if $x(k)<0$;} \\
                \frac{1}{2}x(k)-1, & \ \ \ \ \hbox{if $x(k)>0$.}
              \end{array}
            \right.\]$\ $\\
Clearly $x(k)\in]-1,1[$ for every $k\in\Bbb N$, so system (\ref{eq
example 3 system}) is not ULDI.

Nonetheless system (\ref{eq example 3 system}) is ULI in $1$ step.
Consider indeed the quantization set (defined in Definition \ref{Def
set Q})

\[Q = \bigcup_{i\in\Bbb Z} [i,i+1[\times [i,i+1[\]$\ $\\
and observe that

\[\frac{1}{2}Q+\left(
                 \begin{array}{c}
                   u \\
                   u^{\prime} \\
                 \end{array}
               \right)\ \bigcap\ Q = \emptyset\ \ \forall
               u\not=u^{\prime}.\]$\ $\\
This in turn implies that system (\ref{eq example 3 system}) is ULI
in $1$ step. $\diamondsuit$
\end{Example}

\section{Conclusions}

In this paper we studied left invertibility of output-quantized
linear systems, and we proved that it is equivalent, under suitable
conditions, to left D-invertibility, a stronger notion, much easier
to detect (Theorem \ref{Theor D-inv exp-contr}). More precisely the
condition under which left invertibility and left D-invertibility
are equivalent is that the elements of the dynamic matrix of the
system form an algebraically independent set. Therefore the set of
matrices for which left D-invertibility and left invertibility are
equivalent is a full measure set (Theorem \ref{Theor alg
independent}). Moreover there is a standard way to create matrices
whose elements are an algebraically independent set (Theorem
\ref{Theor Lin-Weier}). Notice that algebraic conditions play a
central role in investigation of left invertibility of quantized
systems as well in other fields when a quantization is introduced
(see for instance \cite{Bic,Chi}).

Future research will include further investigation on the
equivalence between left invertibility and left D-invertibility to
matrices whose elements are not algebraically independent.

\section*{Appendix}

\subsection*{\textbf{Proof of Proposition \ref{Prop transverse bounded}}}

Define the following sequence of sets:

\[\left\{\begin{array}{l}
    S_0 =  Ker(\pi_{d-p})\ \ \ (= \langle e_{p+1},\ldots,e_{d}\rangle)\\
    S_{i+1} = A(S_i)\cap Ker(\pi_{d-p}).
  \end{array}\right.
\]$\ $\\
Then $A$ does not have an invariant subspace included in
$Ker(\pi_{d-p})$ if and only if $S_{d-p} = \{0\}$. To prove this,
first of all observe that $S_i$ is a subspace of $Ker(\pi_{d-p})$
for every $i\in\Bbb N$. Moreover it holds that $S_i$ is a subspace
of $S_{i-1}$ for every $i = 1,\ldots$. To prove this, by induction:
\begin{itemize}
  \item $S_1 = A(S_0)\cap S_0$ is clearly a subspace of $S_0$;
  \item Suppose that $S_i$ is a subspace of $S_{i-1}$. Then $S_{i+1} = A(S_i)\cap
  Ker(\pi_{d-p})$ is a subspace of $A(S_{i-1})\cap
  Ker(\pi_{d-p}) = S_i$.
\end{itemize}

Suppose that there exists a subspace $W\subset Ker(\pi_{d-p})$. Then
$W\subset A^i[Ker(\pi_{d-p})]$ for every $i=0,\ldots,d-p$. This in
turn implies that $S_{d-p}\not=\{0\}.$

Viceversa, suppose that $S_{d-p}\not=\{0\}$. Since the length of the
sequence $\{S_0,\ldots,S_{d-p}\}$ is $d+1$, there exists
$j\in\{0,\ldots,d-p\}$ such that $dim (S_j) = dim(S_{j+1}) =
dim[A(S_j)\cap Ker(\pi_{d-p})]$, where we indicate with $dim$ the
dimension. But $S_{j+1}$ is a subspace of $S_{j}$, and this in turn
implies that $S_j=S_{j+1}$. Moreover both are subspaces of
$Ker(\pi_{d-p})$. So

\[S_{j+1} = A(S_j)\cap Ker(\pi_{d-p}) = S_j.\]$\ $\\
Therefore $A(S_j)=S_j$ because $dim[A(S_j)]\leq dim(S_j)$. So $S_j$
is an invariant subspace of $Ker(\pi_{d-p})$.\\

Consider now the following succession of sets:
\[\left\{\begin{array}{l}
    \tilde S_0 =  (]-1,1[)^p\times Ker(\pi_{d-p})\\
    \tilde S_{i+1} = A(\tilde S_i)\cap (]-1,1[)^p\times Ker(\pi_{d-p}).
  \end{array}\right.
\]$\ $\\
Since $S_{d-p} = \{0\}$ (a point) $\tilde S_{d-p}$ must be bounded.
Finally, observe that $\tilde S_{d-p}$ is exactly the set of
possible states $z(d-p)$ when every $z(0),\ldots,z(d-p)$ is in
$(]-1,1[)^p\times Ker(\pi_p)$. The Proposition is thus proved.
$\diamondsuit$\\

\textbf{Proof of Theorem \ref{Theor alg independent}}:

\begin{Definition}\label{def_I0_Omega}
Let us parametrize the possible pairs of states
$(x,x^{\prime})\in\Bbb R^{2d}$ such that $x^{\prime}-x\in
(]-1,1[)^p\times\langle e_{p+1},\ldots,e_{d}\rangle$ with the set

\[Q^{\prime} = \left\{ \left(t_1,\ldots,t_d,t_1+s_1,\ldots,t_p+s_p,\tilde{t}_{p+1},\ldots,\tilde{t}_{d}\right):
\ t_i,\tilde{t}_j\in\Bbb R,s_k\in]-1,1[ \right\}.\]$\ $\\
Moreover define

\[Q^{\prime}_i = \{(t,t+s)\in\Bbb R^2:\ t\in\Bbb R,\ s\in]-1,1[\}\subset \langle e_i,e_{d+i}\rangle,\]

\[Q_i = \bigcup_{j\in\Bbb Z}\  [j,j+1[\times[j,j+1[\subset \langle
e_i,e_{d+i}\rangle.\]$\ $\\
If
$X=(t_1,\ldots,t_d,t_1+s_1,\ldots,t_p+s_p,\tilde{t}_{p+1},\ldots,\tilde{t}_d)\in
Q^{\prime}$, for $i=1,\ldots,p$ define $d_i(X)$ to be the distance,
measured along the line
\[\{t_1,\ldots,\underbrace{\tau_i}_{varies},\ldots,t_d,t_1+s_1,\ldots,\underbrace{\tau_i}_{varies}+s_i,\ldots,t_p+s_p,\tilde{t}_{p+1},\ldots,\tilde{t}_d:
\ \tau_i\in\Bbb R\}\] from the set $\Omega_i = \{X\in\Bbb R^{2d}:\
\varpi_jX=0\ for\ j\not=i,i+d\}$. $\diamondsuit$
\end{Definition}

\begin{Lemma}\label{Lemma ij kron cond}
Fix a sequence $\{U(j)\}_{j\in\Bbb N}$ of inputs for the doubled
system (\ref{Sistema2d}). Suppose that $\forall\epsilon>0$, $\forall
m\in\Bbb N$, $\forall s_1,\ldots,s_p\in]-1,1[$, there exists
$t_1,\ldots,t_d,\tilde{t}_{p+1},\ldots,\tilde{t}_d\in\Bbb R$ such
that, if $\{X(j)\}_{j=0}^m\subset Q^{\prime}$ is the orbit of the
doubled system (\ref{Sistema2d}) with
$X(0)=(t_1,\ldots,t_d,t_1+s_1,\ldots,t_p+s_p,\tilde{t}_{p+1},\ldots,\tilde{t}_d)$,
then the following holds
\begin{equation}\label{eq-ij kron cond}
frac\left( \frac{d_i(X(j))}{\sqrt{2}} \right)<\epsilon,
\end{equation}
for every $i=1,\ldots,p$, $j=1,\ldots,m$. Then the system is not
ULI.
\end{Lemma}

\emph{Proof:} Suppose that an orbit $\{X(i)\}_{i=1}^{\infty}$ of the
doubled system is included in $Q^{\prime}$ and consider the
2-dimensional plane spanned by $\langle e_i,e_{d+i}\rangle$. Observe
that $frac\left( \frac{d_i(X(j))}{\sqrt{2}} \right)=0$ if and only
if $X(j)$ belongs to some translation of $\Omega_i$ along the
bisecting line of the 2-dimensional plane $\langle
e_i,e_{d+i}\rangle$, that is entirely included in $Q_i$, i.e. a
translation that takes $\Omega_i$ to the ``bottom-left boundary'' of
a square of $Q_i$.

Suppose now that the relation (\ref{eq-ij kron cond}) is satisfied
for every $\epsilon$. It's now easy to see that, for every $X\in
Q^{\prime}$ there exists $\epsilon>0$ such that, if $frac\left(
\frac{d_i(X)}{\sqrt{2}} \right)<\epsilon$ then the projection of $X$
on $Q_i^{\prime}$ is indeed in $Q_i$ (thanks to the periodicity of
the property of ``being in $Q_i$'', see also the illustrative
example). Therefore, if the relations (\ref{eq-ij kron cond}) are
satisfied, then there exists an arbitrary long orbit included in
$Q$. $\diamondsuit$

\begin{Proposition}\label{Prop_ind_matr_pow}
Suppose that the entries of a matrix $A=\{a_{ij}\}\in\Bbb R^{d\times
d}$ are an algebraically independent set, and denote with
$a^{(k)}_{ij}$ the entries of the matrix $A^k$. Then the set
\begin{equation}\label{eq_set_lin_ind}
\left\{ a^{(k)}_{ij}:\ i,j=1,...d;\ k=1,\ldots,N \right\}
\end{equation}
is a linearly independent set for every $N\in\Bbb N$.
\end{Proposition}

\emph{Proof:} First, note that all $a^{(k)}_{ij}$ are polynomials of
degree $k$ in the $a_{ij}$'s. Since the $a_{ij}$ are algebraically
independent, they can be treated formally as the independent
variables of polynomials in $d^2$ variables (more precisely there
exists a ring isomorphism between $\Bbb
Q[a_{11},a_{12},\ldots,a_{dd}]$ and the ring of polynomials in $d^2$
variables $\Bbb Q[x_{1},\ldots,x_{d^2}]$, see \cite{Art}). If a
nontrivial linear combination of the elements of the set
(\ref{eq_set_lin_ind}) is zero, then there exists a nontrivial
polynomial in the $a_{ij}$ which is zero, so there exists a $k_0$
such that a nontrivial linear combination of the $a_{ij}^{(k_0)}$'s,
seen as \emph{polynomials in the $a_{ij}$'s}, which is zero. These
are the entries of the matrix $A^{k_0}$, so there would exist a
nontrivial linear relation among these entries. Suppose this is the
case. If this linear relation results in a linear relation among
polynomials which is not identically zero, we are done. Indeed, if
there exists $\lambda_{ijk_0}\in\Bbb Q$ and
$x_{11},\ldots,x_{dd}\in\Bbb R$ such that (note that the entries
$a_{ij}^{(k_0)}$ of the matrix $A^{k_0}$ are seen as polynomials in
the variables $a_{11},\ldots,a_{dd}$, renamed as
$x_{11},\ldots,x_{dd}$)

\[\sum_{i,j} \lambda_{ijk_0}a_{ij}^{(k_0)}(x_{11},\ldots,x_{dd})\not= 0,\]
then, substituting the $a_{ij}$'s to the $x_{ij}$'s, it is not
possible that

\[\sum_{i,j} \lambda_{ijk_0}a_{ij}^{(k_0)}(a_{11},\ldots,a_{dd})=
0,\] since the $a_{ij}$'s are algebraically independent.

Therefore we only have to show that it is not possible that

\[\sum_{i,j} \lambda_{ijk_0}a_{ij}^{(k_0)}(x_{11},\ldots,x_{dd})\equiv 0,\]
i.e. that this polynomial cannot be identically zero. Now note that
the matrices $M=\{m_{ij}\}\in\Bbb R^{d\times d}$ whose entries do
not satisfy the (nontrivial) linear relation $\sum_{i,j}
\lambda_{ijk_0}m_{ij}= 0$ form a full measure set, dense in $\Bbb
R^{d\times d}$. On the other hand also the matrices with distinct
eigenvalues form a full measure set, dense in $\Bbb R^{d\times d}$.
Therefore there exists a matrix $\overline M=\{\overline m_{ij}\}$
having distinct eigenvalues, whose entries do not satisfy the linear
relation $\sum_{i,j} \lambda_{ijk_0}\overline m_{ij}= 0$.

Since $\overline{M}$ has distinct eigenvalues, there exists a matrix
$B$ such that $B^{k_0}=\overline{M}$ (diagonalize and take
$k_0$-roots of the eigenvalue). Denote with $b^{(k_0)}_{ij}$ the
entries of the matrix $B^{k_0}$: then the $b^{(k_0)}_{ij}$'s do not
satisfy the linear relation $\sum_{i,j}
\lambda_{ijk_0}b_{ij}^{(k_0)}= 0,$ since $B^{k_0}=\overline{M}$.
This implies that $\sum_{i,j}
\lambda_{ijk_0}a_{ij}^{(k_0)}(x_{11},\ldots,x_{dd})$ is not
identically zero as a polynomial. $\diamondsuit$\\

We now prove Theorem \ref{Theor alg independent}. Consider the state
$X(m)$ of the doubled-system (\ref{Sistema2d}) at instant $m$ given
by an initial condition $X(0)\in Q^{\prime}$ and an input sequence

\[U(1)=\left(
        \begin{array}{c}
          u(1) \\
          u^{\prime}(1) \\
        \end{array}
      \right),\ldots,U(m)=\left(
                            \begin{array}{c}
                              u(m) \\
                              u^{\prime}(m) \\
                            \end{array}
                          \right)
      \in\mathcal{U}\times\mathcal{U}:\]
Then

\[X(m) = \left(
           \begin{array}{cc}
             A & 0 \\
             0 & A \\
           \end{array}
         \right)^mX_0+\ldots+\left(\begin{array}{c}
                            Bu(m) \\
                            Bu^{\prime}(m) \\
                      \end{array}\right)\]
Suppose that the ($d-$dimensional) system (\ref{sistema}) is not
uniformly left D-invertible. So there exists arbitrarily long orbits
of the ($2d-$dimensional) doubled-system (\ref{Sistema2d}) included in $Q^{\prime}$.\\
In the following we provide conditions such that
$\forall\epsilon>0$, $\forall m\in\Bbb N$, $\forall
s_1,\ldots,s_p\in]-1,1[$, there exists $t_1,\ldots,t_d,\tilde
t_{p+1},\ldots,\tilde t_d\in\Bbb R$ such that, if $\{X(j)\}_{j=0}^m$
is the orbit of the doubled-system (\ref{Sistema2d}) with
$X(0)=(t_1,\ldots,t_d,t_1+s_1,\ldots,t_p+s_p,\tilde{t}_{p+1},\ldots,\tilde{t}_d)$,
then the following holds
\begin{equation}
frac\left( \frac{d_i(X(j))}{\sqrt{2}} \right)<\epsilon,
\end{equation}
for every $i=1,\ldots,p$, $j=1,\ldots,m$. Therefore the system will
be not uniformly left invertible by Lemma \ref{Lemma ij kron cond}.
These conditions will be verified by a full measure set. Consider
the set
\[\Bbb S^{\prime} = \left\{ A\in\Bbb R^{d\times d}:
\{a_{ij}\}_{i,j=1}^{d}\ is\ an\ algebraically\ independent\
set\right\}.\] Set $A\in\Bbb S^{\prime}$. For $i=1,\ldots,p$,\ \ \
$\varpi_{\langle e_i,e_{d+i}\rangle}X(j)$ has the form

\[\left(
    \begin{array}{c}
      \varpi_iX(j) \\
      \varpi_{i+d}X(j) \\
    \end{array}
  \right) = \pi_{\langle i,i+d\rangle}\left[
            \left(
              \begin{array}{cc}
                A & 0 \\
                0 & A \\
              \end{array}
            \right)^jX(0)+\left(
                            \begin{array}{cc}
                              A & 0 \\
                              0 & A \\
                            \end{array}
                          \right)^{j-1}U(1)+\ldots+U(j)\right] = \]

\[ = \pi_{\langle i,i+d\rangle}\left[
                \begin{array}{c}
                \left(\begin{array}{cccccc}
                        a_{11} & \ldots & a_{1d} & 0 & \ldots & 0 \\
                        \vdots & \ddots & \vdots & \vdots & \ddots & \vdots \\
                        a_{d1} & \ldots & a_{dd} & 0 & \ldots & 0 \\
                        0 & \ldots & 0 & a_{11} & \ldots & a_{1d} \\
                        \vdots & \ddots & \vdots & \vdots & \ddots & \vdots \\
                        0 & \ldots & 0 & a_{d1} & \ldots & a_{dd}
                      \end{array}\right)^j\ \ \left(
                                            \begin{array}{c}
                                            t_1 \\
                                            \vdots \\
                                            t_d \\
                                            t_1+s_1 \\
                                            \vdots \\
                                            t_p+s_p \\
                                            \tilde t_{p+1} \\
                                            \vdots \\
                                            \tilde t_{d} \\
                                            \end{array}
                                        \right)\ +\
                                        const.\\\end{array}\right] = \]

\[ =  \left(
  \begin{array}{l}
     c_{i,1}^{(j)}\ t_1+\ldots+c_{i,d}^{(j)}\ t_d \\
     \\
     c_{i,1}^{(j)}\ t_1+\ldots+c^{(j)}_{i,p}\ t_p+c^{(j)}_{i,p+1}\ \tilde{t}_{p+1}
+\ldots+c_{i,d}^{(j)}\ \tilde{t}_d \\
  \end{array}
\right)\ + const.\] where $c_{i,l}^{(j)}$ is the entry $(i,l)$ of
the matrix $A^j$. By Proposition \ref{Prop_ind_matr_pow} the set
$\left\{ a^{(j)}_{il}:\ i,l=1,...d;\ j=1,...N \right\}$ is a
linearly independent set, so, by Kronecker's Theorem (Theorem
\ref{Theorem Kro 2}) there exists a choice of $(t_1\ldots,t_d,\tilde
t_{p+1},\ldots,\tilde t_d)$ such that equation (\ref{eq-ij kron
cond}) is satisfied, and Lemma \ref{Lemma ij kron cond} thus apply.

To prove that the set of matrices with algebraically independent
entries are a full measure set, first observe that the set of
polynomial $P\in\Bbb Q[\zeta_1,\ldots,\zeta_{d^2}]$ is countable.
For a single polynomial $P$ the set
\[0_P = \left\{ (x_1,\ldots,x_{d^2})\in\Bbb R^{d^2}:\ P(x_1,\ldots,x_{d^2})=0 \right\}\]
is a finite union of manifolds of dimension at most $d^2-1$. So the
measure of $0_P$ is zero. Moreover
\[\Bbb S^{\prime} = \bigcup_{P\in\Bbb Q[\zeta_1,\ldots,\zeta_{d^2}]} 0_P,\]
i.e. $\Bbb S^{\prime}$ is a countable union of sets of measure zero,
which in turn implies that the measure of $\Bbb S^{\prime}$ is zero.
$\diamondsuit$\\

\section{Notations}

The authors made an effort to simplify notations, though they are
intrisically complex. For this reason in this ``special'' section we
collect the notations used in this paper, ordered as their
appearance.

\begin{itemize}
\item $\mathcal P = \bigcup\mathcal{P}_{i}$ is the uniform partition (Definition \ref{Def_unif_partition});

\item $\mathcal U$ is the finite alphabet of inputs (just after the
Definition \ref{Def_unif_partition});

\item $q:(x\in \mathcal P_i) \mapsto i$ is the output quantizer (just after the Definition \ref{Def_unif_partition});

\item $Q$ is quantization set (Definition \ref{Def set
Q});

\item $\mathcal A$ denotes the attractor of a system, (Definition \ref{Def_attractor});

\item $\mathcal I$ denotes the invariant set of a system (Definition
\ref{Def_attractor});

\item $\phi$ indicates the function that associates to each input
sequence its limit point (Theorem \ref{Teorema insieme fisso});

\item $G_k,IG_k$ are respectively the graph of depth $k$ associated to the attractor
$\mathcal A$, and the internal invertibility graph (Definitions
\ref{Def. Graph associaed to attractor}, \ref{Def Int Ext Inv
Graph});

\item $V_{IG_k}$ denotes the union of vertices (which are sets)
of the internal invertibility graph (Definition \ref{Def Int Ext Inv
Graph});

\item $\mathcal V = \mathcal{U}-\mathcal{U}=\{u-u^{\prime}:\
u\in\mathcal U,\ u^{\prime}\in\mathcal U\}$ is the input set of the
difference set (Definition \ref{Def Difference System})

\item $D_{k_1}^{k_2}(z(0),v(1),\ldots,v(k_2))$ denotes the sequence
$(\pi_pz(k_1),\ldots,\pi_pz(k_2))$ generated by the difference
system with initial condition $z(0)$ and inputs $v(1)\ldots,v(k_2)$
(Definition \ref{Def_Dk1k2});

\item $\Bbb S_D(B,\mathcal U)$ denotes the set of matrices $A\in\Bbb R^{d\times
d}$ such that the system (\ref{sistema}) is ULDI (before Theorem
\ref{Theor alg independent});

\item $\Bbb S(B,\mathcal U)$ denote the set of matrices $A\in\Bbb R^{d\times
d}$ such that the system (\ref{sistema}) is ULI (before Theorem
\ref{Theor alg independent});

\item $Q^{\prime} = \left\{
\left(t_1,\ldots,t_d,t_1+s_1,\ldots,t_p+s_p,\tilde{t}_{p+1},\ldots,\tilde{t}_{d}\right)^T:t_i,\tilde{t}_j\in\Bbb
R,s_k\in]-1,1[ \right\}$.

\item $d_i(X)$ is the measure of a distance defined in Definition
\ref{def_I0_Omega};

\item $\Omega_i$ is the union of the two coordinate axes of $\langle
e_i,e_{d+i}\rangle$ (Definition \ref{def_I0_Omega});

\end{itemize}

$\ $\\
$\ $\\
\end{document}